\documentclass[11pt,a4paper]{article}
\usepackage[utf8]{inputenc}
\usepackage{amsmath}
\usepackage{mathrsfs}
\usepackage{amsfonts}
\usepackage{amssymb}
\usepackage[dvips]{graphicx}
\usepackage{listings}
\usepackage{hyperref} 
\usepackage{color}
\usepackage{makerobust}
\usepackage[table]{xcolor}
\usepackage{amsthm}
\usepackage{xspace}
\usepackage{dsfont} % for mathds
\usepackage{subfigure}

\usepackage[T1]{fontenc}
\usepackage{amsfonts} % need for mathbb
\usepackage[geometry]{ifsym} % triangle
\usepackage{savesym} % for savesymbol, restoresymbol command needed to resolve conflict between ifsym and bbding
\savesymbol{Cross}
\savesymbol{Square}
\savesymbol{TriangleUp}
\savesymbol{TriangleDown}
\usepackage{bbding} % ellipse
\restoresymbol{bb}{Cross}
\restoresymbol{bb}{Square}
\restoresymbol{bb}{TriangleUp}
\restoresymbol{bb}{TriangleDown}
\usepackage{dsfont} % for mathds
\usepackage[font=footnotesize,labelfont=bf]{caption} % small captions
\usepackage{authblk}

\usepackage{tikz}
\usetikzlibrary{patterns,snakes,calc,positioning}
\usepackage{verbatim}

\newtheorem{definition}{Definition}[section]
\newtheorem{assumption}{Assumption}[section]

\newcommand{\ve}[1]{
    \mathbf{#1}
}

\newcommand{\eps}{\varepsilon}
\newcommand{\wos}{\texttt{WoS}\xspace}
\newcommand{\mlwos}{\texttt{MLWoS}\xspace}
\newcommand{\Wtot}{W_L^\text{tot}}      % total amount of work done for a simulation with L+1 levels l=0,...,L

\newcommand{\YY}[1]{
    Y_{\eps_{#1}}
}        % wos estimator for u(x) on level l
\newcommand{\EE}[1]{                % expectation value
    \mathds{E}\!\left[#1\right]
}
\newcommand{\EEx}[1]{                % expectation value
    \mathds{E}_x\!\left[#1\right]
}
\newcommand{\EM}[1]{                % approximation of expectation value
    E_M\!\left[#1\right]
}
\newcommand{\EMl}[2]{               % approx. of expectation value on level l
    E_{M_#1}\!\left[#2\right]       % \EMl{l-1}{Y}
}
\newcommand{\var}[1]{               % variance
    \mathds{V}\text{ar}\!\left[#1\right]
}
\newcommand{\PP}[1]{
    \mathds{P}\!\left[#1\right]
}
\newcommand{\LL}[1]{                % L^2(\Omega) norm
    \lVert#1\rVert_{L^2(\Omega)}
}
\newcommand*\semicircle{
    \raisebox{-1.2pt}{
        \begin{tikzpicture}
            \draw (0.12,0) arc (0:180:0.12) -- cycle;
        \end{tikzpicture}
}}

\newcommand{\Dsquare}{D_\text{\Square}}

\newcommand{\Dhemisphere}{D_{\!\!\!\semicircle}\!\!}

\newcommand{\dd}[1]{\,\mathrm{d}#1} % integration: \int f(x) \dd{x}
\newcommand{\BB}{\mathcal{B}}       % Borel sigma-algebra
\newcommand{\RR}{\mathds{R}}        % set of real numbers
\newcommand{\NN}{\mathds{N}}        % set of natural numbers
\newcommand{\OO}{\mathcal{O}}       % big O notation (asymptotic complexity)
\newcommand*{\defeq}{\mathrel{\vcenter{\baselineskip0.5ex \lineskiplimit0pt
                     \hbox{\scriptsize.}\hbox{\scriptsize.}}}%
                     =}

\renewcommand{\a}{\ve{a}}
\renewcommand{\b}{\ve{b}}
\renewcommand{\c}{\ve{c}}

\renewcommand{\r}{\ve{r}}

\newcommand{\x}{\ve{x}}

\def\addsymbol #1: #2#3{$#1$ \> \parbox{115mm}{#2 \dotfill \pageref{#3}}\\}
\def\newnot#1{\label{#1}}

\title{Multilevel Monte Carlo for the Feynman-Kac Formula for the Laplace Equation}
\author[1,2]{Stefan Pauli\thanks{The work of this author has been funded by the ETH
      interdisciplinary research grant CH1-03 10-1.}}
      \author[1]{Robert Gantner}
      \author[1]{ Peter Arbenz}
      \author[3] {Andreas Adelmann}  
\affil[1]{ETH Zurich, Computer Science Department, 8092 Zurich, Switzerland}
\affil[2]{ETH Zurich, Seminar for Applied Mathematics, 8092 Zurich, Switzerland}      
\affil[3]{PSI, 5232 Villigen, Switzerland}   
\begin{document}
\maketitle

\section*{Abstract}
Since its formulation in the late 1940s, the Feynman-Kac formula has proven
to be an effective tool for both theoretical reformulations and practical
simulations of differential equations. The link it establishes between such
equations and stochastic processes can be exploited to develop Monte Carlo
sampling methods that are effective, especially in high dimensions. There exist many techniques of improving standard Monte Carlo sampling methods, a relatively new development being the so-called Multilevel Monte Carlo method. This paper investigates the applicability of multilevel ideas to the stochastic representation of partial
differential equations by the Feynman-Kac formula, using the Walk on Sphere algorithm to generate the required random paths. 
We focus on the Laplace equation, the simplest elliptic PDE, while mentioning some extension possibilities.

\section{Introduction}
The Monte Carlo (MC) error normally converge like $1/\sqrt M$, where $M$ is
the number of samples. If all samples were equally expensive the error
versus work convergence rate would be $1/\sqrt {W}$, where $W$ is the total work. 
Note that for a large class of problems the accuracy of the MC method is not only determined by the number of samples used, but also by the discretization error accepted in the computation of each sample.
A small discretization error, using e.g.~the Walk on Spheres (\wos) algorithm~\cite{Muller}, comes with larger computation costs per sample. 
Hence the error
versus work convergence rate falls short of $1/\sqrt {W}$. 
If applicable, Multilevel Monte Carlo (MLMC) methods may reach this
optimal $1/\sqrt {W}$ convergence rate.
This motivated us to evaluate the Feynman-Kac formula with MLMC methods instead of MC methods.

In this paper we present a procedure to evaluate a Feynman-Kac formula with MLMC using the Walk on Spheres method. As a model problem  we use the Laplace equation, 
which is solved in high dimensions for instance for option pricing \cite[Chap. 8]{hilber2013}
 or in particle accelerator simulations~\cite{Adelmann2010}.
We prove that the error versus work converges with the optimal $1/\sqrt {W}$ convergence rate, and compute the suboptimal convergence when using MC methods. The  MLMC method outperforms the MC method only by a $\log(W)$ term. Our MLMC simulations, executed with an MPI parallel implementation, where up to twice as fast compared to the standard MC implementation.

First, in Section~\ref{sec:WoS} we introduce the Feynman-Kac formula, using the \wos algorithm, and compute the error versus work convergence rate when using MC. Then, we derive a method to use MLMC in this setting, and evaluate its error versus work convergence rate in Section~\ref{sec:multilevel_walk_on_spheres}. 
In Section~\ref{sec:numResults} we present numerical results, which quantify the advantage of using MLMC. 
Finally, we draw our conclusions in Section~\ref{cha:outlook_and_conclusion}.

\section{Standard Walk on Spheres} \label{sec:WoS}
%%               %%
% Brownian Motion %
%%               %%
\subsection{Brownian motion} % (fold)
A Brownian motion~\cite{Chang2007} denoted by $X_t$, started at a point in a certain connected domain $D$, has several characteristic quantities. These include the first exit time and the first exit point, two important concepts in the application of the Feynman-Kac formula. 
\begin{definition}[First exit time]\label{def:exit_time}
The time at which a realization of a Brownian motion $X_t$, started at some point $x\in D\subset\RR^d$, first hits the domain boundary $\partial D$ is called the first exit time $\tau \defeq \inf \{t>0: X_t \in \partial D\}$.
\end{definition}
\begin{definition}[First exit point]\label{def:exit_point}
    The point at which a realization of a Brownian motion $X_t$, started at some point $x\in D\subset\RR^d$, first hits the domain boundary $\partial D$ is called the first exit point $X_\tau$.
\end{definition}

Of great importance for the derivation of the Walk on Spheres algorithm is the \textit{distribution} of first exit points of a Brownian motion in a ball $B_{x}(R)$\newnot{symbol:Bxr} of radius $R$ centered at the point $x$.
If the starting point of the Brownian motion is the center of the ball, this distribution is uniform on its spherical surface, independent of the dimension
\cite[Theorem~3]{Wendel1980}.

In evaluating functionals of stochastic processes, one is often interested in their expectations and variances.
Throughout this paper, the $L^2$-norm over the sample space $\Omega$ is defined as the the expectation $\LL{X}:= \EE{|X|^2}^{\frac12}$.
% section brownian_motion (end)

%%             %%
% Elliptic PDEs %
%%             %%

\subsection{Feynman-Kac formula for the Laplace equation}

The Feynman-Kac formula, developed by Richard Feynman and Mark Kac \cite{Kac1949}, gives a probabilistic representation of the solution to certain PDEs at a single fixed evaluation point $x$. 
The Feynman-Kac formula allows us to write a very general elliptic PDE
\begin{align*}
    Lu &= \frac12 \Delta u + \sum_{i=1}^d b_i(x) \partial_i u + c(x)u = -g(x), & &\text{ in }D\subset\RR^d, \\
    u(x) &= f(x), & &\text{ on } \partial D,
\end{align*}
along with the assumption that $b$, $c$ and $g$ are smooth and satisfy a Lipshitz growth condition \cite{Buchmann2003} as an expectation over a stochastic process $X_t$ beginning at the point $x$.
This process fulfills the stochastic differential equation $\dd{X}_t=b(x)\dd{t} + \!\dd{W}$ where $\dd{W}$ is a Brownian increment, and is stopped as soon as it hits the boundary.
Denoting the first exit time by $\tau$ and the condition that $X_0=x$ by the subscript $x$ in $\mathds{E}_x$, the formula can be written as
\begin{equation}\label{eq:feynman_kac}
    u(x) = \EEx{ f(X_\tau)\exp\left(\int_0^\tau\!c(X_s)\dd{s}\right) + \int_0^\tau\!g(X_t)\exp\left(\int_0^t\!c(X_s)\dd{s}\right)\dd{t} }.
\end{equation}
In this paper we are only interested in the Laplace equation $\Delta u = 0$ in $d$ dimensions on a domain $D\subset\RR^d$ with boundary values $u(x) = f(x)$ on $\partial D$. This simplifies the development of an MLMC method, which can be used as a basis for more general equations. 
In this simple case the Feynman-Kac representation is given as an expectation over the exit point $X_\tau$ of a Brownian motion started at the point $x$. Denoting the condition $X_0=x$ by the subscript $x$ in $\mathds{E}_x$, the solution is written as~\cite{Buchmann2003}
\begin{equation}
	u(x)=\EEx{f(X_\tau)}=:\EEx{Y}.
	\label{eq:feynman_kac}
\end{equation}

\subsection{Walk on Spheres (\wos)} % (fold)
\label{subsec:WoS}
The Walk on Spheres (\wos) algorithm \cite{Muller} can be viewed as an alternative to a conventional detailed simulation of the drift-free Brownian motion inside the domain $D$. Started at the center of a sphere the exit point distribution is known. This allows to simulate the Brownian motion using discrete jumps of a given size. Starting at $\tau=0$ with $X_0=x$ the algorithm measures the distance from the current position $X_t$ to $\partial D$ and jumps this distance in a uniformly random direction to the next position $X_{t+1}$. The algorithm terminates when it is $\eps$ close to the boundary $\partial D$ at the point $X_N$, where $N$ denotes the number of \wos steps needed. The first exit point is approximated by the point $\overline{X}_N \in \partial D$ that is closest to $X_N$.
This \wos algorithm may be used to solve the Laplace equation~\cite{Muller}.
We recall that $\EE{Y}=u(x)$ when $Y$ is computed with a realization of a processes starting at $x$.
$Y_\eps\defeq f(\overline{X}_N),$ is the estimator based on an \wos processes with the discretization parameter $\eps$. 
The discretization error presented in \cite{Mascagni2003} is given as
\begin{equation}\label{eq:err_discretization}
    e_\text{discr} = \LL{\EE{Y}-\EE{Y_\eps}} \le \OO(\eps).
\end{equation}

In \cite{Binder2009}, a worst-case upper bound of $\OO(\eps^{4/d-2})$ for three or more dimensional problems ($d\ge3$) is derived for the expected \wos path length $\EE{N}$. For domains fulfilling certain regularity conditions tighter upper bounds of $\EE{N} = \OO(\log^p(\eps^{-1}))$, $d\ge2$ are proven for $p=1$ or $p=2$ depending on the domain. These estimates are relevant for bounding the expected work per sample.

\subsection{Monte Carlo}

\subsubsection{Statistical error} % (fold)
\label{ssub:statistical_error}
Given $M$ realizations $\{Y_\eps^i\}_{i=1}^M$ of the random variable $Y_\eps$ obtained from independent paths of the \wos process, the value for $\EE{Y_\eps}$ can be approximated by the estimator $\EM{Y_\eps} = \frac1M\sum_{i=1}^MY_\eps^i$.
The error of this Monte Carlo estimator can be written in terms of the random variable $Y_\eps$ as
\begin{align}\label{eq:err_statistical}
    \begin{split}
        e_\text{stat}   &= \LL{\EE{Y_\eps}-\EM{Y_\eps}} = \var{\frac1M\sum_{i=1}^M Y_\eps^i}^\frac12 \\
                        &= \frac1{\sqrt{M}}\sqrt{\var{Y_\eps}} \le \frac1{\sqrt{M}}\LL{Y_\eps},
    \end{split}
\end{align}
where we use the fact that the samples $Y_\eps^i$ are independent and identically distributed (i.i.d.)\ realizations of $Y_\eps$.

\subsubsection{Total error} % (fold)
\label{ssub:total_error}
The total error of an estimation $\EM{Y_\eps}$ of $\EE{Y}$ by the Walk on Spheres algorithm can be written using the triangle inequality of the $L^2(\Omega)$ norm as
\begin{equation*}
    %\EE{\left|\EE{Y} - \EM{Y_\eps}\right|^2}^{\frac12} \le \EE{\left|\EE{Y}-\EE{Y_\eps}\right|^2}^\frac12 + \EE{\left|\EE{Y_\eps}-\EM{Y_\eps}\right|^2}^\frac12,
    \LL{\EE{Y} - \EM{Y_\eps}} \le \LL{\EE{Y}-\EE{Y_\eps}} + \LL{\EE{Y_\eps}-\EM{Y_\eps}},
\end{equation*}
resulting in a total error of
\begin{equation}\label{eq:err_total}
    e_\text{tot} \defeq \LL{\EE{Y} - \EM{Y_\eps}} \le \OO(\eps) + \sqrt{\frac{\var{Y_\eps}}{M}}.
\end{equation}
The two terms at the right follow by \eqref{eq:err_discretization} and \eqref{eq:err_statistical}, respectively.
% subsubsection total_error (end)

\subsubsection{Error equilibration} % (fold)
\label{ssub:error_equilibration}
For a fixed prescribed error, one way of choosing the sample size is by equilibrating the statistical and discretization errors in~\eqref{eq:err_total}, $\OO(\eps) = \sqrt{\var{Y_\eps}/M}$.
This yields the relationship
\begin{equation}\label{eq:M_equilibrium}
    M=\OO(\eps^{-2}),
\end{equation}
giving a total error behavior of $\OO(\eps)$.
% subsubsection error_equilibration (end)

\subsection{Error vs. expected work} % (fold)
\label{sub:error_vs_work}
The expected total work $\EE{W}$ of a \wos simulation is the number of paths times the expected length of a path:
\begin{equation}
    \EE{W} = M\cdot \EE{N}.
\end{equation}

Each bound for the expected path length $\EE{N}$ shown in section~\ref{subsec:WoS} is multiplied with $M$ derived in \eqref{eq:M_equilibrium}, and the resulting expected work is solved for $\eps$. Equation \eqref{eq:Wvseps_WOS_a} is valid for $d\ge3$, equation \eqref{eq:Wvseps_WOS_b} for well-behaved domains $d\ge2$.  To derive \eqref{eq:Wvseps_WOS_b} we must make use of the Lambert $W$-function, $\mathcal{W}_\text{lam}(\cdot)$, defined as the inverse  of the map $w\mapsto w\exp(w)$. It can be approximated with the truncated expansion $\mathcal{W}_\text{lam}(x) \approx \log(x)-\log\log(x)$ \cite{Bruijn1970}.
This yields the following relationships:
\begin{subequations}
    \begin{align}
        \EE{W} &= \OO(\eps^{4/d-4}) &   e_\text{tot} = \OO(\eps) &= \OO( \EE{W}^{\frac14\frac{d}{1-d}})\label{eq:Wvseps_WOS_a}\\
        \EE{W} &= \OO(\eps^{-2}\log^2(\eps^{-1})) & e_\text{tot} = \OO(\eps) &= \OO(\EE{W}^{-\frac12}\log(\EE{W})) \label{eq:Wvseps_WOS_b}
    \end{align}
\end{subequations}
The total error is linear in $\eps$ by the choice in \eqref{eq:M_equilibrium}. In both cases the \wos algorithm performs worse than the optimum possible in a Monte Carlo setting, namely $W^{-\frac12}$. This motivates the formulation of a multilevel version of the \wos algorithm, with the hopes of achieving the optimal convergence rate.

\section{Multilevel Walk on Spheres (\mlwos)} % (fold)
\label{sec:multilevel_walk_on_spheres}
In this section, a multilevel version of the \wos algorithm is formulated and its error behavior analyzed.
The main idea is to execute the \wos algorithm on different ``discretization levels'', meaning for different values of the discretization parameter $\eps$.
The subscript $\ell\in\{0,\ldots,L\}$ is used to denote a certain discretization level, where $\ell=0$ is the coarsest discretization level, corresponding to a \wos simulation with $\eps_0$, and $\ell=L$ is the finest discretization level.

\subsection{Multilevel formulation} % (fold)
\label{sub:multilevel_formulation}
On each discretization level $\ell=0,\ldots,L$, we define the discretization parameter $\eps_\ell\defeq\eta^{-\ell}\eps_0$, \newnot{symbol:epsl}\newnot{symbol:eta} where $\eta>1$, i.e.~the width of the stopping region is divided by $\eta$ between successive levels. 

In the following definition, a single \wos process with $\eps_\ell$ is used to define the multilevel process in order to incorporate the fact that stopping points on higher levels are continuations of previously stopped processes.
\begin{definition}[Multilevel Walk on Spheres process]\label{def:mlwos}
    Given a domain $D\subset\RR^d$, a point $x\in D$, a discretization parameter $\eps_0$ and an $\eta>1$.
    Consider a \wos process $\{X_i\}_{i=0}^{N_\ell}$ started at $x$ with $\eps_\ell = \eta^{-\ell}\eps_0$ and the pair $\left( X_{N_{\ell-1}}, X_{N_\ell} \right)$ obtained by $N_{\ell-1} = \min \{N\in\NN : d_{\partial D}(X_N)<\eps_{\ell-1}\}$, where $d_{\partial D}(X_N)=\min_{x'\in\partial D}|X_N-x'|$ is the distance to the boundary $\partial D$. 
    The Multilevel Walk on Spheres process (\mlwos) on the level $\ell$ is the set of all such pairs.
    %Then, the Multilevel Walk on Spheres process (\mlwos) is defined as the totality of such sets.
\end{definition}

In the context of multilevel Monte Carlo, we have multiple estimators $\YY{\ell}$, one for each discretization level with discretization parameter $\eps_\ell$.
These are obtained in the same way as for the non-multilevel case, e.g.~$\YY{\ell}=f(\overline{X}_{\!N_\ell})$. The expectation of the estimator $\YY{L}$ is written in a multilevel form as 
\begin{equation}\label{eq:mlmc_exact}
    \EE{\YY{L}} = \EE{\YY{0}} + \sum_{\ell=1}^L \EE{\YY{\ell}-\YY{\ell-1}}.
\end{equation}
Replacing the expectation with the average over $M_\ell$ realizations on each discretization level $\ell$, we get the MLMC estimate
\begin{equation}\label{eq:mlmc_approx}
    E\left[\YY{L}\right] = \EMl{0}{\YY{0}} + \sum_{\ell=1}^L \EMl{\ell}{\YY{\ell}-\YY{\ell-1}}.
\end{equation}
%$E[\cdot]$ (without a subscript) denotes this concept of using $M_\ell$ samples on level $\ell$.
The ``sample'' on level $\ell>0$ is now $(\YY{\ell}-\YY{\ell-1})$, for which it is assumed that the two values on discretization levels $\ell$ and $\ell-1$ come from the same \wos path.
The implication is that they should be simulated with the \mlwos process from \autoref{def:mlwos}.
The remaining estimator on discretization level $0$ is computed with an ordinary \wos simulation.

Note that the individual $\YY{\ell}$ are estimators on the \textit{discretization level} $\ell$, whereas $\YY{\ell}-\YY{\ell-1}$ are referred to as estimators on \textit{level} $\ell$.
% subsection multilevel_formulation (end)

\subsection{Multilevel error bounds} % (fold)
\label{sub:multilevel_error_bounds}
The multilevel Monte Carlo error is defined as the difference between the expectation of the exact estimator $\EE{Y}$ and the MLMC approximation involving all levels, $E\left[\YY{L}\right]$, and is given by~\cite{Barth2011}
\begin{align*}
    \LL{\EE{Y}-E[\YY{L}]} &\le \LL{\EE{Y}-\EE{\YY{L}}} \\
    & \qquad + M_0^{-\frac12}\LL{\YY{0}} + \sum_{\ell=1}^L M_\ell^{-\frac12}\LL{ \YY{\ell}-\YY{\ell-1} }.
\end{align*}

\subsection{Asymptotic variance convergence rate} % (fold)
\label{sub:variance_convergence_rate}
We want to bound the variance of the estimator on a level $\ell$ in the multilevel error bound. 
In more exact terms, we want to determine for which functions $f$ and domains $D$ the relationship
\begin{equation}\label{eq:varconvrate}
    \LL{\YY{\ell}-\YY{\ell-1}} = \OO(\eps_\ell^s) \quad\forall \ell>\ell_{min}
\end{equation}
holds for some $\ell_{min}>0,s>0$.
This behavior is desired since it allows a good multilevel performance as the variance of the estimators on fine levels is already small, thus requiring less fine-level realizations.

Conditioning the expectation on the current position $X_{N_{\ell-1}}$ incorporates the fact that the simulation continues with a given path from level $\ell-1$, not generating a completely independent one on the finer level $\ell$. Using $\overline{X}_{N_\ell} = \overline{X}_{N_{\ell-1}}+\Delta x$ yields
\begin{equation}\label{eq:varconv_1}
    \LL{\YY{\ell}-\YY{\ell-1}}^2 = \EE{\left.|f(\overline{X}_{N_{\ell-1}}+\Delta x)-f(\overline{X}_{N_{\ell-1}})|^2\right|X_{N_{\ell-1}}}.
\end{equation}
In order to write this expectation in terms of $|\Delta x|$, not in terms of a function of the process we make the following assumption:
%Differentiability of $f$ would allow this to to be rewritten using the mean-value theorem, however a more general assumption is given using Hölder continuity.
\begin{assumption}[Hölder continuity]
    There exist $C,\alpha > 0$ such that
    \begin{equation*}
        |f(x)-f(y)| \le C|x-y|^\alpha
    \end{equation*}
    for all $x,y\in\partial D$.
\end{assumption}\newnot{symbol:alpha}
This assumption implies
\begin{align}\label{eq:varconv_2}
    \EE{\left.|f(x+\Delta x)-f(x)|^2\right|x} &\le C^2\EE{\left.|\Delta x|^{2\alpha}\right|x}.
\end{align}
The conditioning can be omitted by bounding the right-hand side of \eqref{eq:varconv_2} by the maximum over $x$, the stopping point of the process on level $\ell-1$.

Splitting the expectation into a converging part where $|\Delta x| \le R_\ell$ and a ``diverging'' part with $R_\ell<|\Delta x| \le |D|$, where $|D|$ is the diameter of the circumsphere of $D$, we get \newnot{symbol:absD}
\newcommand{\pdiv}{p_{\text{div}}}
\begin{align}\label{eq:Edx_ansatz}
    \begin{split}
    \EE{|\Delta x|^{2\alpha}} &\le \PP{|\Delta x_\ell| \le R_\ell}R_\ell^{2\alpha} + \PP{|\Delta x_\ell| > R_\ell}|D|^{2\alpha} \\
                              &= (1-\pdiv)R_\ell^{2\alpha} + \pdiv |D|^{2\alpha}.
    \end{split}
\end{align}
Since the $|D|$ in the second term is a domain-dependent constant, we must find a bound for the divergence probability $p_{\text{div}}$ in terms of the discretization parameter $\eps$.
As shown below, this is possible with the resulting behavior $\pdiv\propto\frac{\eps_\ell}{R_\ell}$.

In order for the expectation to converge with a certain rate in $\eps$, both terms $R_\ell^{2\alpha}$ and $\pdiv$ should be of the same order.
Equating the two, solving this equation for $R_\ell$ and inserting back into \eqref{eq:Edx_ansatz} yields
\begin{align}
    R_\ell &= C\eps_\ell^{\frac1{2\alpha+1}} \nonumber \\
    \EE{|\Delta x|^{2\alpha}} &\le C\eps_\ell^{\frac{2\alpha}{2\alpha+1}} \label{eq:varconv_result}, \\ 
    \LL{\YY{\ell}-\YY{\ell-1}}^2&\le C\eps_\ell^{\frac{2\alpha}{2\alpha+1}}. \nonumber
\end{align}
For differentiable $f$, $\alpha=1$ this yields the convergence rate $\eps_\ell^\frac13$ in \eqref{eq:varconvrate}.

\subsubsection{Convergence of divergence probability} % (fold)
\label{ssub:convergence_of_divergence_probability}
Our goal here is to obtain a bound on the divergence probability $\pdiv$ on a certain level $\ell$ depending on $\eps_\ell$.
Divergence here means that the process does not result in an estimate $\overline{X}_{N_{\ell+1}}$ that is located within a ball $B_{\overline{X}_{N_\ell}}\!(R_\ell)$ of radius $R_\ell$ around the current projected stopping point $\overline{X}_{N_{\ell}}$.

In the derivation we use the so-called harmonic measure.
There exist a few equivalent definitions of harmonic measure which are relevant to the current application.
The first defines the harmonic measure as a harmonic function satisfying certain boundary conditions,
the second is instructive in the context of Monte Carlo approximations of an integral and
the third uses a definition relying on the distribution of the first exit point given a Brownian motion. Consider the closed set $D\subset\RR^d$ with $d\ge2$ and let $\BB(\partial D)$ denote the $\sigma$-algebra of subsets of $\partial D$ and define:

\begin{definition}[Harmonic measure -- Dirichlet solution]\label{def:harmonicmeasure_dirichlet}
    The harmonic measure $\omega_D: D\times\BB(\partial D) \to [0,1]$, viewed as a function $\omega_D(x,E)$ of $x$ for every fixed $D$ and $E$, is the unique harmonic function that satisfies the boundary condition \cite[p. 39]{Garnett2005}
    \begin{equation*}
        f(x) =
            \begin{cases}
                1, &\text{ if }x\in E,\\
                0, &\text{ if }x\in \partial D \backslash E.
            \end{cases}
    \end{equation*}
    %In this case, uniqueness must be shown with a stronger version of the maximum principle shown by Lindelöf
\end{definition}

\begin{definition}[Harmonic measure -- integral representation]\label{def:harmonicmeasure_integral}
    The harmonic measure $\omega_D: D\times\BB(\partial D) \to [0,1]$ is the unique function that satisfies \cite[sec. 4.3]{Ransford1995}
    \begin{item}
        \item[(a)] for each $x\in D$, $E\mapsto\omega_D(x,E)$ is a probability measure on $\partial D$
        \item[(b)] if $f:\partial D\to\RR$ is continuous, the solution of the Laplace equation in $D$ is given by
            \begin{equation}\label{eq:u_harmonicmeasure}
                u(x) = \int_{\partial D}\! f(x') \,\omega_D(x,\!\dd{x'}).
            \end{equation}
    \end{item}
    %In other words, the contributions of the boundary condition to the solution are weighted according to the harmonic measure.
\end{definition}
\begin{definition}[Harmonic measure -- hitting probability]\label{def:harmonicmeasure_hitting}
Let $X_t^x$ denotes a Brownian motion started at $x$.
    The harmonic measure $\omega_D: D\times\BB(\partial D) \to [0,1]$ is given by
    \[ \omega_D(x,E) \defeq \PP{X_\tau^x \in E, \tau<\infty}, \]
    where $\tau = \inf\{t\ge0:X_t^x\in\partial D\}$ is the first exit time of $X_t^x$ from $D$.
    \cite[Theorem F.6.]{Garnett2005}
\end{definition}

Writing $\omega_D(x,x') = \frac{\partial {\omega}_D(x,x')}{\partial {x'}}$\newnot{symbol:omega_density} as a probability density function of the distribution of first exit points $x'$ of $\partial D$ by a Brownian motion started at $x$, \eqref{eq:u_harmonicmeasure} corresponds to the expectation value of $f(X_\tau)$.
This agrees with the Feynman-Kac formula \eqref{eq:feynman_kac} applied to the Laplace equation.
Simulating realizations of the stochastic process $X_t$ thus amounts to approximating the integral \eqref{eq:u_harmonicmeasure} by Monte Carlo integration, i.e. generating realizations of the exit points distributed according to the harmonic measure.

Often, the harmonic measure in a certain domain can't be given analytically.
Thus, it is important to be able to bound it from above and below.
This can be accomplished using the general inequalities popularized by Carleman, illustrated in \autoref{fig:carleman} and often referred to as \textit{Carleman's Principle of Monotonicity} or the \textit{Principle of Extension of Domain}.
\begin{definition}[Carleman's principle]\label{def:carlemans_principle}
    The harmonic measure of a portion $\alpha$ of the boundary of $G\subset\RR^d$ increases if the region is increased by modifying $\beta=\partial G \backslash \alpha$:
    \begin{equation}\label{eq:carleman_lower}
        \omega_G(x,\alpha) \le \omega_{G'}(x,\alpha).
    \end{equation}
    Since $\omega_G(x,\alpha)+\omega_G(x,\beta)=1$, it follows that the harmonic measure on the extended portion of the boundary decreases with respect to the new domain:
    \begin{equation}\label{eq:carleman_upper}
        \omega_G(x,\beta) \ge \omega_{G'}(x,\beta').
    \end{equation}
    %This is the formulation from \cite[p. 68]{nevanlinna1960analytic}; \cite[p. 131]{doob2001classical} states it for subsets of $\RR^d$.
    \cite[p. 68]{nevanlinna1960analytic}, \cite[p. 131]{doob2001classical}
    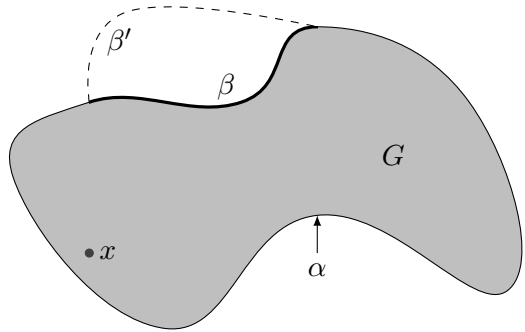
\begin{figure}[ht]
\centering

\begin{tikzpicture}
    %% fill (do this first)
    \begin{scope}
        \clip plot[smooth cycle,tension=0.8] coordinates{(0,1)  (2,1)   (3,2)     (5,1)      (5.5,-1.5)  (3,-0.5) (1,-2) (-1,0)};
        \fill[black!25]
            (current bounding box.south west)
            rectangle
            (current bounding box.north east);
    \end{scope}
    %% Domains (normal and eps)
    \draw[dashed] plot[smooth,tension=0.8] coordinates{(0,1) (0.5,2.2) (3,2)};
    \draw plot[smooth cycle,tension=0.8] coordinates{(0,1)  (2,1)   (3,2)     (5,1)      (5.5,-1.5)  (3,-0.5) (1,-2) (-1,0)};
    \begin{scope}
        \clip (0,0.6) rectangle (3,2);
        \draw[very thick] plot[smooth cycle,tension=0.8] coordinates{(0,1)  (2,1)   (3,2)     (5,1)      (5.5,-1.5)  (3,-0.5) (1,-2) (-1,0)};
    \end{scope}
    %\draw[dashed] plot[smooth cycle,tension=0.7] coordinates{(0.4,0.6) (2.8,1.5) (4.5,0.75) (5.1,-1.6)  (1.2,-2) (-0.55,-0.85)} -- cycle;
    \draw (0.4,1.8) node{$\beta'$};
    \draw (1.8,1.2) node{$\beta$};
    %\draw (3,-1.25) node{$\alpha$};
    \draw[latex-] plot (3,-0.75) -- (3,-1) node[anchor=north]{$\alpha$};
    \draw (4,0.3) node{$G$};
    \filldraw [black!75] (0,-1) circle (1.5pt);
    \draw (0,-1) node[anchor=west]{$x$};
\end{tikzpicture}

\caption{An illustration of Carleman's Principle in two dimensions. $\alpha$ is the unchanged part of the boundary, $\beta$ is changed to $\beta'$.
The extended domain is $G'$, $G$ is the original one.}
\label{fig:carleman}
\end{figure}
.
    %\todo{assume $G$ greenian? (see doob, p. 27, 13. Greenian Sets)}
\end{definition}

Writing the divergence probability $\pdiv$ in terms of a harmonic measure  $\pdiv = \PP{\overline{X}_{N_{\ell+1}}\notin B_{\overline{X}_{N_\ell}}\!(R_\ell)} = \omega_D(X_{N_\ell},\partial D \cap B_{\overline{X}_{N_\ell}}\!(R_\ell)^\mathsf{c})$, we can appeal to \textit{Carleman's Principle}.
This allows the transformation of a general bounded convex domain onto a semi-ball, the harmonic measure of which is known analytically.
Simple trigonometric observations then result in a bound depending on $\eps_\ell$.

\begin{figure}[ht]
\centering

\begin{tikzpicture}
    %% constants
    \def\r{20mm}
    \def\a{17}
    %\coordinate (a) at (-1,-1); % midpoint of semicircle
    %% clip for D'
    \begin{scope}
        \clip ($(-1,-1)+({\r*sin(\a)},-{\r*cos(\a)})$) arc (-90+\a:90+\a:\r);
        \clip plot[smooth cycle,tension=0.8] coordinates{(0,1)     (3,2)     (5,1)   (5.5,-2) (1,-2.5) (-1,-1)} -- cycle;
        \fill[black!25]
            (current bounding box.south west)
            rectangle
            (current bounding box.north east);
    \end{scope}
 
    %% brownian path
    \draw[black!50] (-0.55,-0.85) -- (0,0) -- (1.5,0.2) -- (3.5,-0.5) -- (3.2,0.7) -- (4,1.7);
    \filldraw [black!75] (4,1.7) circle (1.5pt) node[anchor=north west,xshift=-3mm]{$X_{\!N_{\ell+1}}$};
    %% Domains (normal and eps)
    \draw         plot[smooth cycle,tension=0.8] coordinates{(0,1)     (3,2)     (5,1)      (5.5,-2)    (1,-2.5) (-1,-1)} -- cycle;
    \draw[dashed] plot[smooth cycle,tension=0.7] coordinates{(0.4,0.6) (2.8,1.5) (4.5,0.75) (5.1,-1.6)  (1.2,-2) (-0.55,-0.85)} -- cycle;
    \draw (3,-1.5) node{$D$};
    \draw (0.5,-1.4) node{$D'$};
    %% stopping points
    \filldraw [black!75] (-1,-1) circle (2pt);% node[anchor=east,yshift=-3mm,xshift=2mm]{$\overline{X}_{N_\ell}$};
    \filldraw [black!75] (-0.55,-0.85) circle (2pt);
    \draw (-0.55,-0.85) node[anchor=west]{$X_{\!N_\ell}$};
    %% Semicircle
    \draw (-1,-1)+({\r*sin(\a)},-{\r*cos(\a)}) arc (-90+\a:90+\a:\r);
    \draw (-1,-1)+({\r*sin(\a)},-{\r*cos(\a)}) -- (-1,-1);
    \draw (-1,-1)+(-{\r*sin(\a)},{\r*cos(\a)}) -- (-1,-1);
    \draw (-1,0.5) node{$S$};

    %% brace: R_\ell
    \draw[
    decoration={
        brace,
        mirror,
        raise=0.2cm
    },
    decorate
    ] (-1,-1)+(-{\r*sin(\a)},{\r*cos(\a)}) -- (-1,-1)
    node [anchor=east, yshift=7.5mm, xshift=-5.5mm] {$R_\ell$};

    %% arrows: eps_\ell
    \draw[thick,latex-]  (5,1)--(5.5,1.25) node [anchor=west]{$\eps_\ell$};
    \draw[thick,-latex]  (4,0.5)--(4.5,0.75);
   
\end{tikzpicture}

\caption{Domain $D$ for simplicity shown in 2 dimensions, with semicircle $S$ of radius $R_\ell$ used to bound the harmonic measure of a point $X_{N_\ell}$ located at most $\eps_\ell$ from the boundary.
The shaded region $D'$ is the intersection of the domain $D$ with the semicircle $S$.
The path to $X_{N_{l+1}}$ is ``divergent'' with respect to $X_{\!N_\ell}$.}
\label{fig:varconv}
\end{figure}
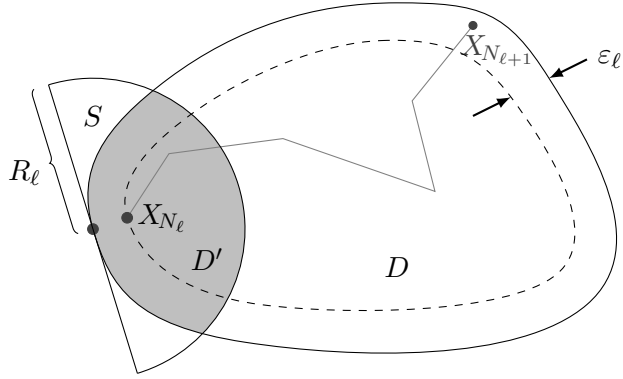

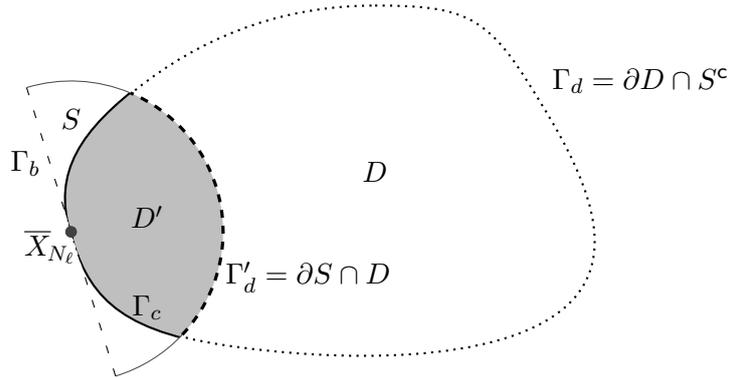
\begin{figure}[ht]
\centering

\begin{tikzpicture}
    %% constants
    \def\r{20mm}
    \def\a{17}
    %% Gamma_d'
    %\draw[loosely dashed] plot[smooth cycle,tension=0.8] coordinates{(0,1)     (3,2)     (5,1)      (5.5,-2)    (1,-2.5) (-1,-1)} -- cycle;
    \draw[dotted,thick] plot[smooth cycle,tension=0.8] coordinates{(0,1)     (3,2)     (5,1)      (5.5,-2)    (1,-2.5) (-1,-1)} -- cycle;
    \draw (5,1) node[anchor=west,xshift=2mm]{$\Gamma_d=\partial D \cap S^\mathsf{c}$};
    %% clip for D'
    \begin{scope}
        \clip (-1,-1)+({\r*sin(\a)},-{\r*cos(\a)}) arc (-90+\a:90+\a:\r);
        \clip[preaction={draw,ultra thick}] plot[smooth cycle,tension=0.8] coordinates{(0,1)     (3,2)     (5,1)   (5.5,-2) (1,-2.5) (-1,-1)} -- cycle;
        \fill[black!25]
            (current bounding box.south west)
            rectangle
            (current bounding box.north east);
    \end{scope}
    %% D,D' label (after intersection) 
    \draw (0,-0.8) node{$D'$};
    \draw (3,-0.2) node{$D$};
    \draw (0,-2) node{$\Gamma_c$};
    %% Gamma_d
    \def\b{30}
    \def\c{40}
    \draw[dashed,very thick] (-1,-1)+({\r*sin(\a+\b)},-{\r*cos(\a+\b)}) arc (-90+\a+\b:90+\a-\c:\r);
    \draw (0.9,-1.6) node[anchor=west]{$\Gamma_d'=\partial S \cap D$};
    %% Semicircle
    \draw[black!70] (-1,-1)+({\r*sin(\a)},-{\r*cos(\a)}) arc (-90+\a:-90+\a+\b:\r);
    \draw[black!70] (-1,-1)+(-{\r*sin(\a)},{\r*cos(\a)}) arc (90+\a:90+\a-\c:\r);
    \draw[loosely dashed] (-1,-1)+({\r*sin(\a)},-{\r*cos(\a)}) -- (-1,-1);
    \draw[loosely dashed] (-1,-1)+(-{\r*sin(\a)},{\r*cos(\a)}) -- (-1,-1);
    \draw (-1,0.5) node{$S$};
    \draw (-1.6,-0.09) node {$\Gamma_b$};
    %\draw (-1,1.3) node {$\partial S\backslash \Gamma_d'$};
    \filldraw [black!75] (-1,-1) circle (2pt);% node[anchor=east,yshift=-3mm,xshift=2mm]{$\overline{X}_{N_\ell}$};
    \draw (-1,-1) node[anchor=east,yshift=-2mm,xshift=2mm]{$\overline{X}_{\!N_\ell}$};
\end{tikzpicture}

\caption{The boundary of the domain $D$ is split into two disjoint portions, the ``convergent'' boundary $\Gamma_c$ and the ``divergent'' boundary $\Gamma_d$.
The portion of the semicircle boundary $\partial S$ inside of $D$ is denoted by $\Gamma_d'$ and the base of the semicircle is $\Gamma_b$.}
\label{fig:varconv_boundary}
\end{figure}

The situation in two dimensions is sketched in \autoref{fig:varconv}: the intersection of the domain $D$ with a semi-ball $S$ of radius $R_\ell>\eps_\ell$ is called $D'$.
In \autoref{fig:varconv_boundary}, the boundaries are split up such that Carleman's method can be applied to the harmonic measure of $\Gamma_d=\partial D \cap B_{\overline{X}_{N_\ell}}(R_\ell)^\mathsf{c}$, the portion of the boundary hit in the case that the process diverges to a point more than $R_\ell$ away from $\overline{X}_{N_\ell}$.

%%In what follows, $x$ is some fixed point in $D'$ and we regard the harmonic measure $\omega_D(x,\cdot)$ as a function of parts of the boundary.
We start by bounding the harmonic measure of $\Gamma_d$ from above by ``shrinking'' $D$ to $D'$ and using the second bound \eqref{eq:carleman_upper}, where $\omega_D(x,\Gamma_d)$ corresponds to the larger domain and thus has smaller harmonic measure:
\begin{equation*}
    \omega_D(x,\Gamma_d) \le \omega_{D'}(x,\Gamma_d').
\end{equation*}
Now, we keep the portion $\Gamma_d'$ fixed and extend the remaining portion $\Gamma_c$ of the boundary of $D'$ to the semicircle boundary $\partial D\backslash\Gamma_d'$.
This lets us use the first inequality \eqref{eq:carleman_lower}, where $\Gamma_d'$ corresponds to $\alpha$ and $\Gamma_c$ corresponds to $\beta$ in the definition, yielding
\begin{equation*}
    \omega_{D'}(x,\Gamma_d') \le \omega_S(x,\Gamma_d') \le \omega_S(x,\partial S\backslash \Gamma_b),
\end{equation*}
where the second inequality comes from adding $\omega_S(x,\partial S\backslash(\Gamma_d'\cup\Gamma_b))$, the harmonic measure of the remaining portion of the arc.

\begin{figure}[ht]
\centering

\begin{tikzpicture}[scale=2.5]
    %% coordinates
    \coordinate (a) at (-1,0);
    \coordinate (b) at (1,0);
    \coordinate (x) at (0,0.3);
    %% semicircle
    \draw [line width=1.5pt] (b) arc (0:180:1);
    \draw (a) -- (b);
    %% angle
    %\def\c{0.35}
    %\coordinate (z) at ($(x)+\c*(0,-1)$);
    %\draw[-latex,black!70] (z) arc (270:330:\c);
    %% labels
    %\draw ($(x)+(0.1,-0.05)$) node[anchor=north]{$\frac\theta2$};
    \filldraw (a) circle (0.2mm) (b) circle (0.2mm) (x) circle (0.2mm) node[anchor=east]{$X_{N_\ell}$};
    \draw ($0.5*(b)+0.5*(a)$) node[anchor=north,xshift=-0.3]{$\overline{X}_{N_\ell}$};
    \filldraw ($0.5*(b)+0.5*(a)$) circle(0.2mm);
    \draw (0.2, 0.8) node{$S$};
    \draw[dashed] (x) -- ($0.5*(b)+0.5*(a)$);
    %\draw (0,0.1) node[anchor=east] {$d$};
    \draw (0.5,0) node[anchor=north] {$R_\ell$};
    %% epsilon
    \draw[dashed,black!60] (-1.2,0.4) -- (1.2,0.4);
    \draw[dotted,black!60] (-1.2,0.2) -- (1.2,0.2);
    \draw[latex-latex] (-0.7,0) -- (-0.7,0.4);
    \draw (-0.7,0.3) node[anchor=west] {$\eps_\ell$};
    %% arrows
    \draw[latex-] (0.6,-0.02) .. controls (0.9,-0.15) .. (1.2,0) node[anchor=west]{$\omega_S(x,\partial S\backslash \Gamma_b)=0$};
    \draw[latex-] (0.707,0.707) -- (1.2,0.707) node[anchor=west]{$\omega_S(x,\partial S\backslash \Gamma_b)=1$};
\end{tikzpicture}
\caption{A semicircle containing the point $X_{N_\ell}$, for simplicity shown for the 2D case. I assumes the value 1 on the arc and 0 on $\Gamma_b$.}
\label{fig_semicircle_construction}
\end{figure}
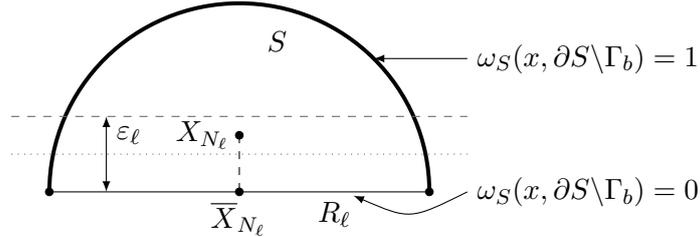

We now assume that $R_\ell \gg \eps_\ell$ and that $x$ lies on the normal of the base at its midpoint. Using that the shortest distance between $x$ and $\Gamma_b$ is less than $\eps_\ell$ allows to use a series approximation of $\omega_S(x,\partial S\backslash \Gamma_b)$ around the midpoint of its base, by using the  Dirichlet solution representation of the harmonic measure Definition~\ref{def:harmonicmeasure_dirichlet}. 
The boundary $\Gamma_b$ is zero by definition, as illustrated in \autoref{fig_semicircle_construction}, hence the same is valid for the first and second derivative in all directions of the base. As by definition the Laplacian is zero, the second derivative in normal direction of the base has to be zero as well. Hence the first derivative in this direction is a positive constant. This results in the final bound for the divergence probability in terms of $\eps_\ell$ 
$$\pdiv \le C \frac{\eps_\ell}{R_\ell},$$
for some constant $C$. 

\subsubsection{Optimal number of samples} 
To simplify notation, we define the estimator $Y_\ell$ to be
\begin{equation*}
    Y_\ell  \defeq \begin{cases}
        \YY{0}, &\text{ if } \ell=0, \\
        \YY{\ell}-\YY{\ell-1}, & \text{ if } \ell>0,
    \end{cases}
\end{equation*}
and $w_\ell$ the expected work performed to compute a sample on level $\ell$.

Following \cite{Muller2012}, the optimization problem to solve is then
\begin{equation}\label{eq:Ml_optimizationproblem}
    \min  \EE{W_\text{tot}}  = \min \sum_{\ell=0}^L M_\ell w_\ell \quad \text{s.t.}\quad \eps_L = \sqrt{\sum_{\ell=0}^L \frac{\var{Y_\ell}}{M_\ell} },
\end{equation}
hence minimizing the expected work, while maintaining a sampling error proportional to the discretization error $(\OO(\eps_L))$. 
\cite{Muller2012} provides the optimal number of samples per level
\begin{equation}
\label{eq:optimal_ML}
    M_\ell = \frac1{\eps_L^2}\sqrt{\frac{\var{Y_\ell}}{w_\ell}}\sum_{\ell=0}^L\sqrt{\var{Y_\ell}w_\ell}.
\end{equation}

\subsection{Error vs. work}
The following basic relations will be used in the sequel.
\begin{subequations}
\begin{align}
    &\text{Refinement of }\eps_\ell     & & \eps_\ell = \eta^{-\ell}\eps_0=\eta^{L-\ell}\eps_L \\
    &\text{Scaling of work}             & & w_\ell = \begin{cases}
       \OO(\eps_\ell^{-\gamma}) \\
       \OO(\log^p(\eps_\ell^{-1})) = \OO(\ell^{{}^p})
    \end{cases} \\
  	&\text{Scaling of variance}             & &  \var{Y_\ell} = \OO(\eps_\ell^{2s})   
\end{align}
\end{subequations}

Determining $M_\ell$~\eqref{eq:optimal_ML} based on the estimates of variance and work lead to 
\begin{align}
 M_\ell &= \begin{cases}
       \eps_L^{-2} \sqrt{\eps_\ell^{2s+\gamma}} \sum_{\ell=0}^L \sqrt{\eps_\ell^{2s-\gamma}}  \\
       \eps_L^{-2} \sqrt{\eps_\ell^{2s}\cdot \ell^{-p}} \sum_{\ell=0}^L \sqrt{\eps_\ell^{2s}\cdot \ell^p}
    \end{cases}      \nonumber \\
&= \begin{cases}
       \eps_L^{-2} \sqrt{\eps_\ell^{2s+\gamma}} \sqrt{\eps_0^{2s-\gamma}} \underbrace{\sum_{\ell=0}^L \sqrt{\eta^{(\gamma-2s)\ell}}}_{\le C, \; \forall 2s>\gamma} \\
       \eps_L^{-2} \sqrt{\eps_\ell^{2s}\cdot \ell^{-p}} \sqrt{\eps_0^2s}\underbrace{\sum_{\ell=0}^L \sqrt{\eta^{-2s\ell}\cdot \ell^p}}_{\le C}
    \end{cases}   \label{eq:Mlopt}
\end{align}
The summands of the remaining sum can be bounded by a convergent geometric sequence. Thus, by the comparison test the infinite series ($L\to\infty$) converges to a constant, in the first case provided $2s>\gamma$.

Inserting \eqref{eq:Mlopt} into the formula $\Wtot = \sum_{\ell=0}^L M_\ell w_\ell$, we can write the work as a function of $\eps_L$, which scales linearly with $e_\text{tot}$
\begin{align*}
    \Wtot &\le  \begin{cases}
       C \eps_L^{-2} \sqrt{\eps_0^{2s-\gamma}} \sum_{\ell=0}^L \sqrt{\eps_\ell^{2s-\gamma}}  \overset{\eqref{eq:Mlopt}}{\le} C^2\eps_L^{-2} \eps_0^{2s-\gamma}, \quad \forall 2s>\gamma\\
       C \eps_L^{-2} \sqrt{\eps_0^{2s}} \sum_{\ell=0}^L \sqrt{\eps_\ell^{2s}\cdot \ell^{p}} \overset{\eqref{eq:Mlopt}}{\le} C^2 \eps_L^{-2} \eps_0^{2s}
    \end{cases} \\
    &= \OO(\eps_L^{-2}).
\end{align*}
Thus, we obtain the optimal convergence rate
\begin{equation}\label{eq:work_vs_error}
    e_\text{tot} \propto (\Wtot)^{-\frac12},
\end{equation}
in the case $w_\ell = \OO(\eps_\ell^{-\gamma})$ provided $2s>\gamma$.

\subsection{Measured values}
\label{sec:Measured_values}
The second approach to determining $M_\ell$ is based on estimating the work $w_\ell$ and the variance $\var{Y_\ell}$ on the different levels, rather than using their asymptotic convergence rates as in~\eqref{eq:Mlopt}. In the absence of alternatives we still use the asymptotic convergence rate of the discretization error $\OO(\eps_L))$. As in~\cite{Muller} the required estimates are computed using the same samples already involved in the MLMC estimator. A certain minimum number of samples (so called warm-up samples) are needed on every level to provide an accurate estimate. Performance disadvantages arise if the number of required warm-up samples exceeds the optimal number of samples~\eqref{eq:optimal_ML}. This typically happens  for the finest level $L$. However, ~\eqref{eq:Mlopt} proves that asymptotically the optimal number of samples grow, which implies that for small $\eps_L$ the required warm-up samples do not exceed the optimal number of samples. This technique complicates the implementation slightly, but does not increase the computational work for this application.

\section{Numerical results}\label{sec:numResults}

\subsection{Model problems}
In this section we formulate model problems, each posed on a different domain to test our algorithm.

\subsubsection{Square} % (fold)
\label{sub:domain_rectangle}
Many of the simplest examples are formulated on  square domains. Here, we use the two-dimensional square domain  $\Dsquare \defeq [0,2]^2$.
In order to more easily measure the convergence rate, a boundary condition resulting in a large discretization error is imposed with Hölder exponent $\alpha =1$. 
\begin{equation}\label{eq:rectangle_BC}
    f(x,y) =
    \begin{cases}
        4(x-\frac12)^2, & \text{if } x \le \frac12, \; 0\le y \le 2,\\
        4(x-\frac32)^2, & \text{if } x \ge \frac32, \; 0\le y \le 2,\\
                     0, & \text{otherwise} , \; 0\le y \le 2.
    \end{cases}
\end{equation}
The midpoint $m=(1,1)$ of the square was chosen as the starting point.

% section domain_rectangle (end)

\subsubsection{Hemisphere} % (fold)
\begin{figure}[ht]
    \centering
    \begin{tikzpicture}[scale=0.8]
        \fill[top color=gray!50!black,bottom color=gray!10,middle color=gray,shading=axis,opacity=0.25] (0,0) circle (4cm and 1cm);
        \fill[left color=gray!50!black,right color=gray!50!black,middle color=gray!50,shading=axis,opacity=0.25] (4,0) arc (0:180:4cm and 4cm) arc (180:360:4cm and 1cm);
        \draw (-4,0) arc (180:360:4cm and 1cm) arc (0:180:4cm and 4cm) -- cycle;
        \draw[densely dashed] (-4,0) arc (180:0:4cm and 1cm);
        \draw (2,0) node[anchor=south] {$R = 1$};
        \draw[-latex] (0,0) -- (4,0);
    \end{tikzpicture}
    \caption[.]{Hemisphere domain $\Dhemisphere$ in three dimensions.}
    \label{fig:domain_hemisphere}
\end{figure}
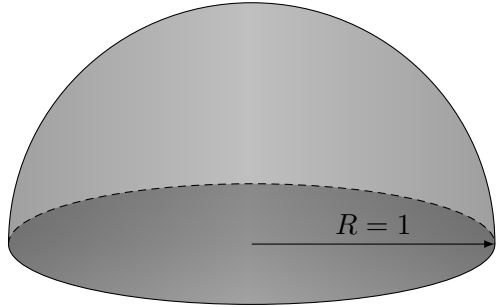

\label{sub:hemisphere}
A nice three-dimensional problem to consider is the Laplace equation on a hemisphere $\Dhemisphere\defeq \{\x\subset\RR^3 \left|\;|\x|\le1, x_3\ge0\right.\}$,\newnot{symbol:Dhemisphere} as illustrated in \autoref{fig:domain_hemisphere}.
The boundary conditions of the Laplace equation are chosen such that the analytical solution is given by $u(\x) = \left[x_1^2+x_2^2+(x_3+1)^2\right]^{-\frac12}$, as stated in \cite{Mascagni2003}. This leads to the following formulation with Hölder exponent $\alpha =1$:
\begin{align}
    \begin{split}
        \Delta u &= 0 \text{ in } \Dhemisphere \\
           u(\x) &= \left[2(x_3+1)\right]^{-\frac12} \text{ for } |\x| = 1 \\
           u(\x) &= [x_1^2 + x_2^2 + 1]^{-\frac12} \text{ for } x_3 = 0.
    \end{split}
\end{align}
The initial point was chosen to be $x_0 \defeq (0.2,0.3,0.1)$.

% section hemisphere (end)

% chapter model_problems (end)

\subsection{Implementation}
A generic C++ implementation for both the conventional and the multilevel version of the Walk on Spheres algorithm was created. The $M_\ell$ are either determined analytically or based on measurements during the simulation. The code uses MPI to generate samples in parallel, such that larger problems can be solved. The code was tested with up to 128 cores. Further information on the implementation is found in~\cite{Gantner2013}.

\subsection{Measurement methodology} % (fold)
\label{sec:Measurement Methodology}
Monte Carlo methods are based on the approximation of an expectation by a sample mean.
The resulting estimator is itself a random variable -- thus, the total error given by
$$\LL{\EE{Y}-E[\YY{L}]}  = \EE{ \| \EE{Y}-E[\YY{L}] \|^2 }^{\frac12}.$$
This outer expectation can again be approximated by a sample mean, corresponding to a repeated call of the corresponding algorithm.
Thus, for a certain set of parameters one must call the implementation with different random number seeds and compute the sample mean over the resulting realizations.
An estimate for a confidence interval can also be computed using the results of these repeated calls.

\subsection{Measured variance convergence rate} % (fold)
\begin{figure}[ht]
    \centering
        \includegraphics[width=0.47\textwidth]{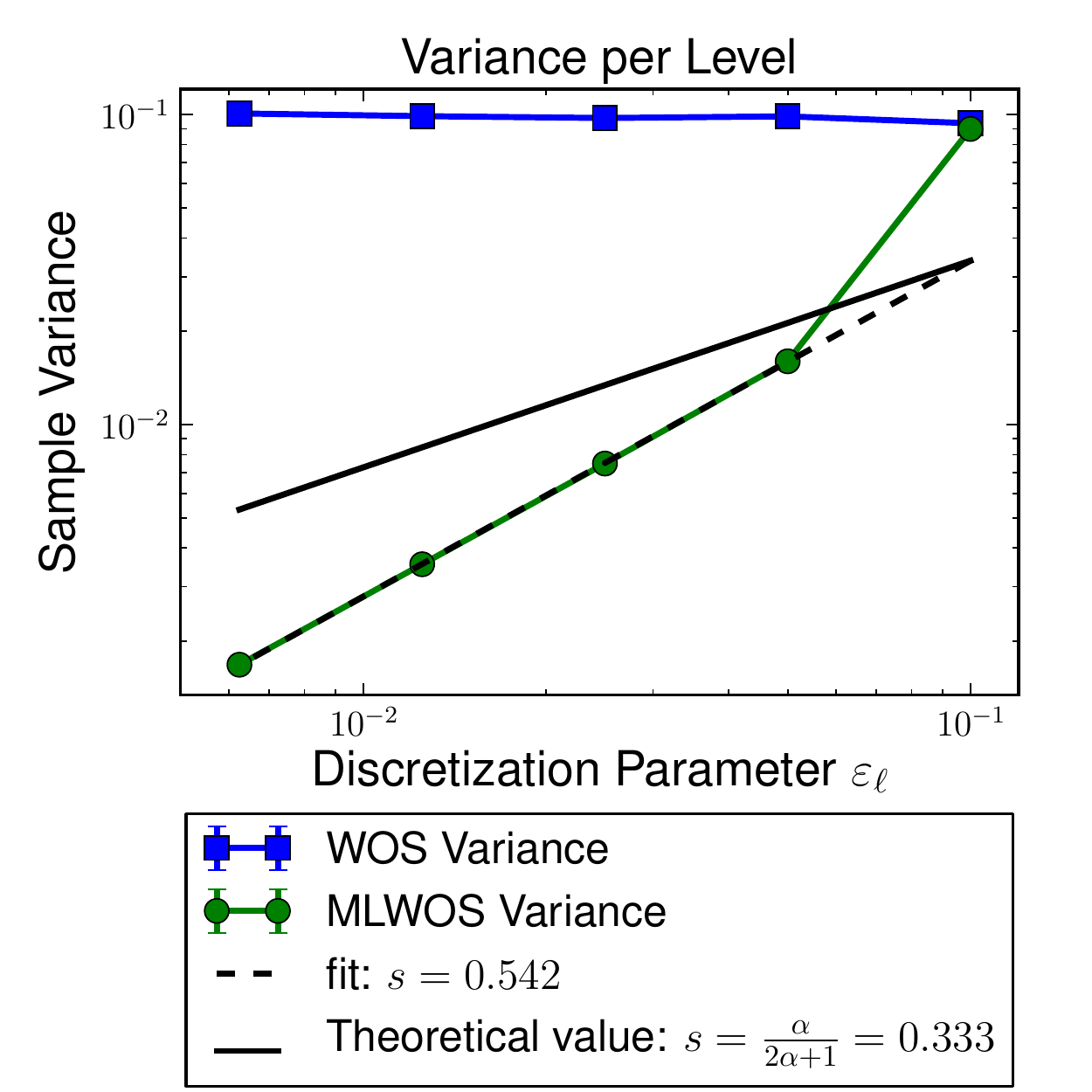}
    \caption[.]{ $\Dsquare$   Convergence of the variance as the discretization parameter $\eps_\ell$ is reduced for the  square domain $\Dsquare$ with $\alpha =1$ using $\eta=2$.}
    \label{fig:res:variance_convergence}
\end{figure}

\label{sec:meas:variance_convergence_rate}
As stated in \eqref{eq:varconvrate}, the variance of the estimate on a level $\ell$ is assumed to converge like $\OO(\eps_\ell^s)$.
In \eqref{eq:varconv_result}, it was shown that $s$ is bounded by $\frac{\alpha}{2\alpha+1}$, $\alpha$ being the exponent of Hölder continuity of the boundary condition.
The results presented here are empirical measurements showing the convergence of the variance of the estimator.
The measurements were averaged over ten calls for the square domain $\Dsquare$.
In \autoref{fig:res:variance_convergence} the theoretical and the measured exponent $s$ are shown.
Shown only for the square domain, but valid for all described domains, the variance convergence rate is around $0.5$, which is better than the value $\frac13$ obtained from the analytical derivation assuming continuity of the boundary condition.
This is not a contradiction since the derivation yields an upper bound.
However, using the analytical upper bound will very likely result in suboptimal results, since obviously one can get away with fewer samples on the fine levels.

% section variance_convergence_rate (end)

\subsection{Measured convergence rate} % (fold)
\label{sec:meas:convergence_rate}

In this section we compare the error versus work of the proposed multilevel \wos algorithm with the plain \wos algorithm. For the multilevel \wos algorithm we distinguish between two strategies, one with analytically derived $M_\ell$, $w_\ell$ and $\var{Y_\ell}$ (MLWOS) and the other with $M_\ell$ based on measured  $w_\ell$ and $\var{Y_\ell}$ (MEAS) according to Section~\ref{sec:Measured_values}.

\begin{figure}[ht]
    \centering
    \subfigure[][$\eta=2$]{
        \includegraphics[width=0.47\textwidth]{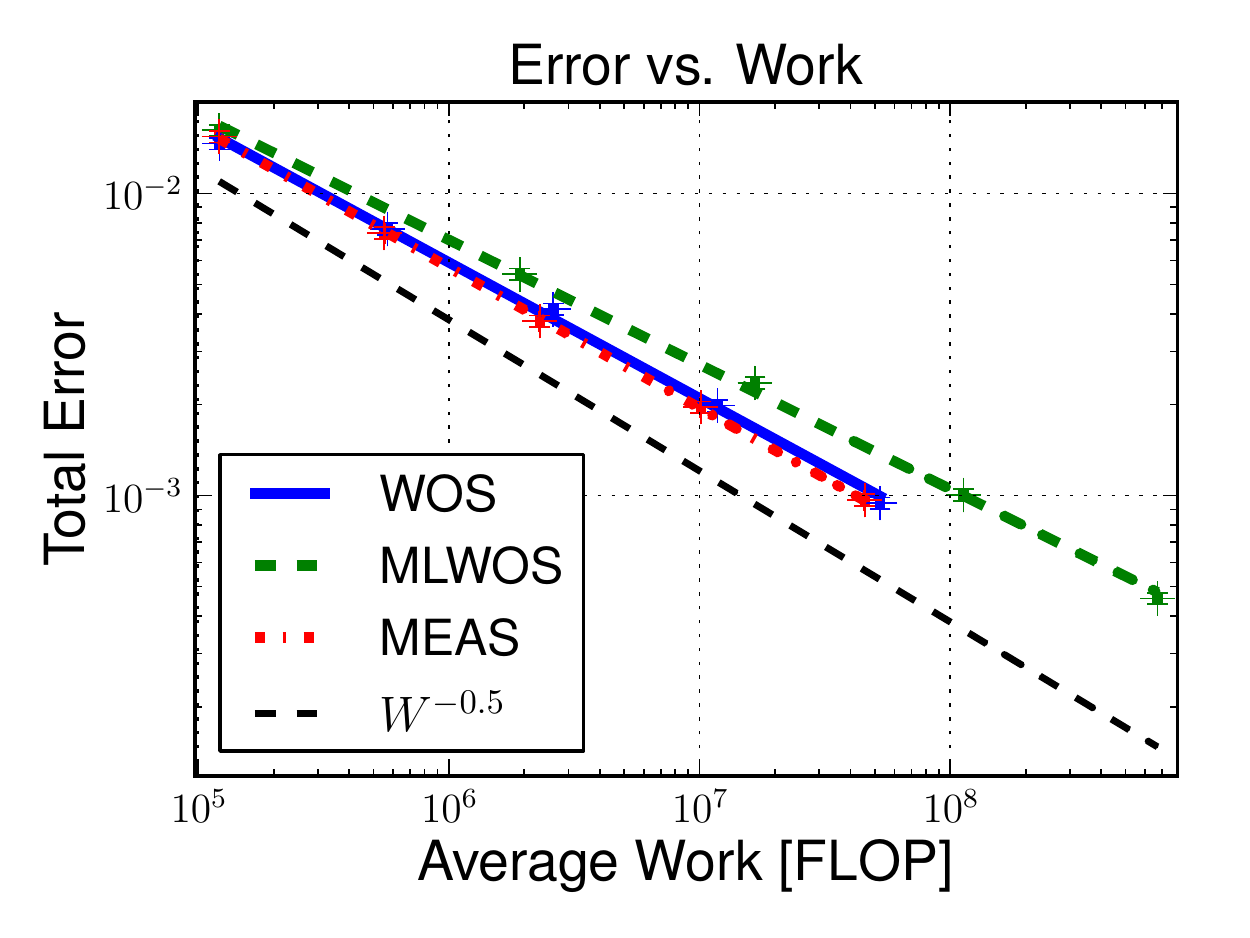}
    } 
    \subfigure[][$\eta=8$]{
        \includegraphics[width=0.47\textwidth]{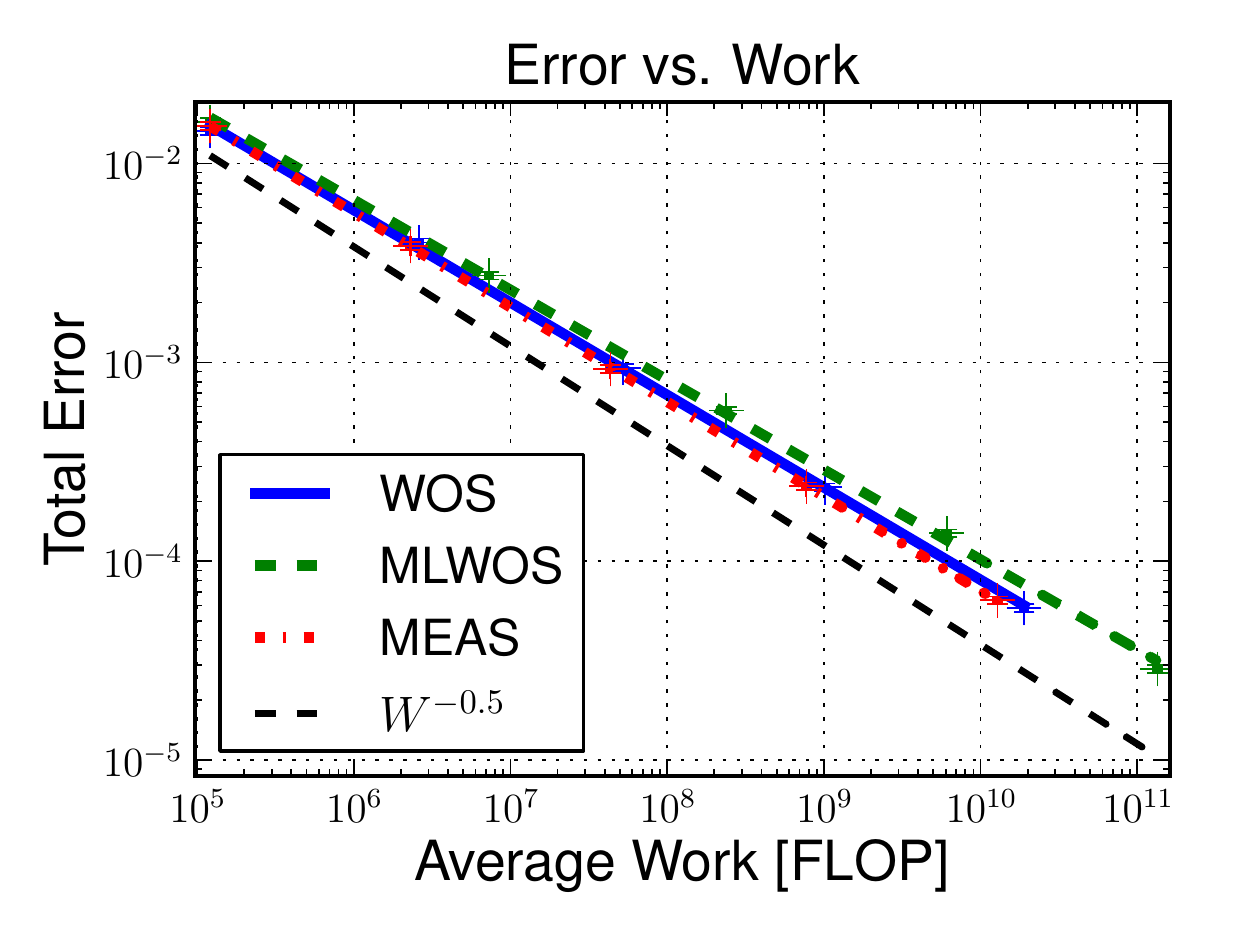}
    }
    \subfigure[][$\eta=16$]{
        \includegraphics[width=0.47\textwidth]{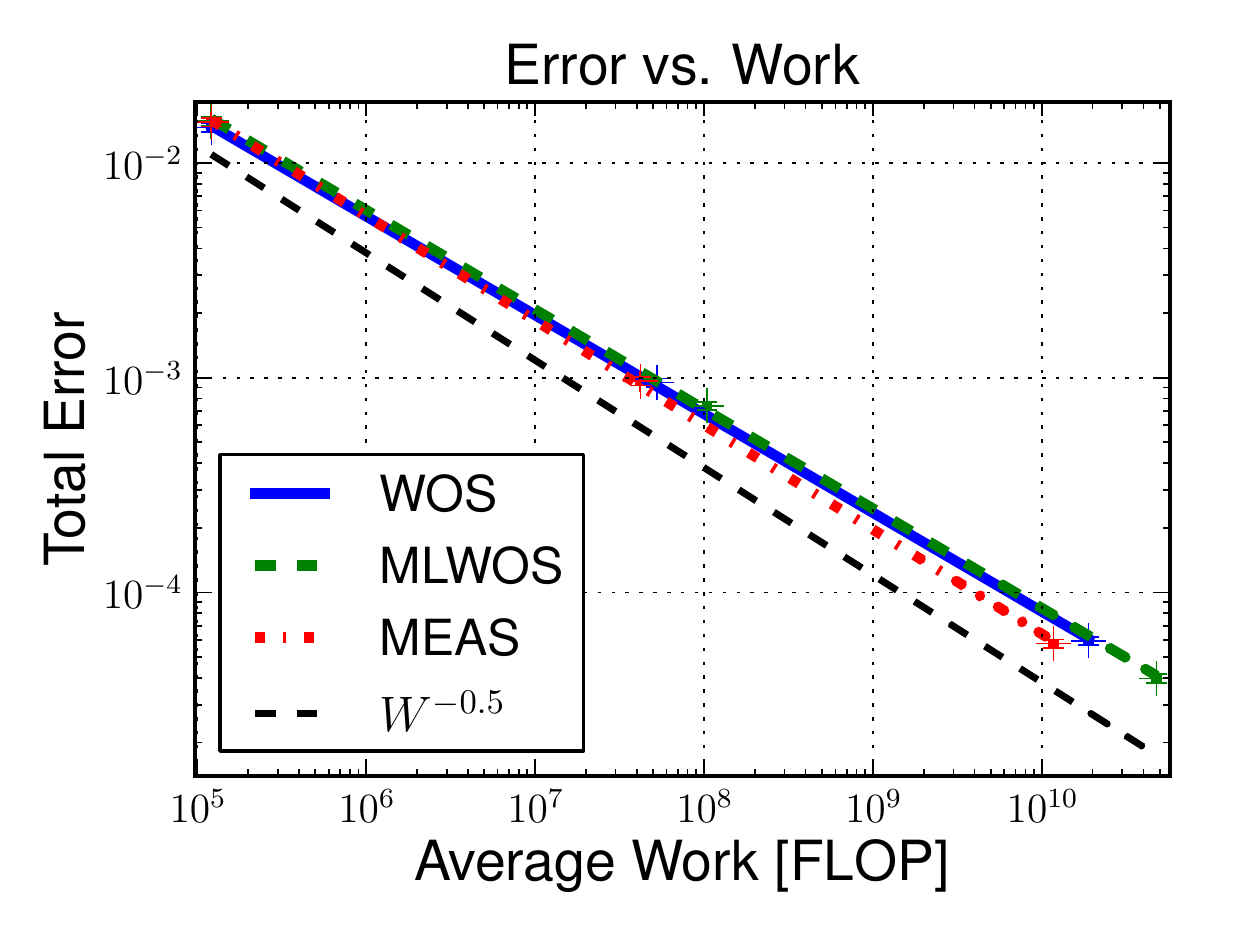}
    }
    \caption[.]{Convergence of error vs.~work for the problem posed on $\Dsquare$, for various values of the refinement parameter $\eta$.}
    \label{fig:res:convergence_eta}
\end{figure}

Measurements performed on the domain $\Dsquare$ are shown in \autoref{fig:res:convergence_eta}. We performed measurements with $\eta=2$, $\eta=8$ and $\eta=16$. The measured values for the average error and the average work are shown, together with the $1\sigma$ confidence interval, for different algorithms. Sampling was performed until the confidence interval was very small. In all cases the multilevel \wos with analytically derived $M_\ell$ (MLWOS) performs poorly compared to the plain \wos algorithm. The MLWOS performance improves for $\eta=16$, but an improvement over the plain \wos algorithm is not measured. The multilevel \wos with $M_\ell$ using measured values (MEAS) always performs better than the plain \wos algorithm, especially for $\eta=16$ where the computation is up to 2 times faster. 

The analytically derived $M_\ell$ are based on the $\var{Y_\ell}$.  As observed in Section~\ref{sec:meas:variance_convergence_rate} $\var{Y_\ell}$ converges faster than predicted in our theory. Hence it is not surprising that the multilevel \wos algorithm with analytically derived $M_\ell$ performs suboptimally. The multilevel \wos algorithm with $M_\ell$ based on measurements does not suffer from this problem, therefore we observe a better performance.  

The parameter $\eta$ does not influence the convergence rate of the multilevel scheme, but a clever choice may asymptotically reduce the work by a constant. Therefore it is expected that the performance of the multilevel \wos algorithm depends on $\eta$. For certain problems (see e.g.~\cite{Giles2008}), this parameter is optimized analytically. Here, various values are tried on the domain $\Dsquare$ in order to empirically find good values. We observed that $\eta=16$ is a better choice than the other tested $\eta=8$ and $\eta=2$.

% various model problems

\begin{figure}[ht]
    \centering
    \includegraphics[width=0.47\textwidth]{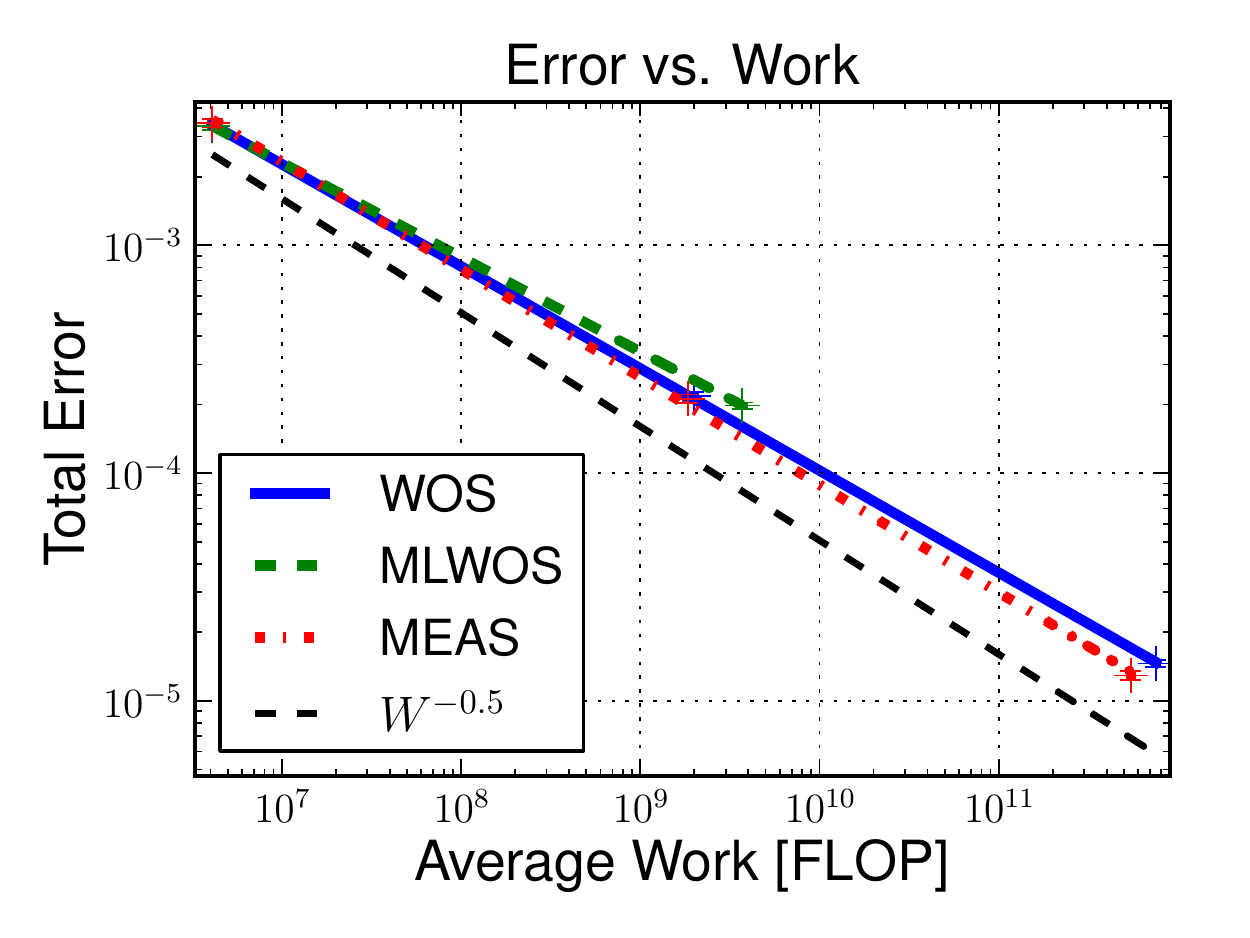}
    \caption[.]{ Convergence of the error vs.~work for the problem posed on $\Dhemisphere$, using $\eta=16$.}
    \label{fig:res:convergence_4}
\end{figure}

\autoref{fig:res:convergence_4} shows measurements for the domain $\Dhemisphere$ for $\eta=16$. The results are similar to the ones seen for the domain $\Dsquare$. 
The multilevel \wos with analytically derived $M_\ell$ (MLWOS) performs poorly, where as the multilevel \wos with $M_\ell$ using measured values (MEAS) is up to 2 times better than the plain \wos algorithm.

% section convergence_rate (end)

% chapter results (end)

\section{Conclusions and outlook} % (fold)
\label{cha:outlook_and_conclusion}
In this work, a first application of the Multilevel Monte Carlo method to the Feynman-Kac formula was developed.
A novel version of the Walk on Spheres process, the ``Multilevel Walk on Spheres'' was formulated and central quantities were derived.
It was proven that the ``Multilevel Walk on Spheres'' algorithm convergence rate of error versus work is optimal. 

In order to measure the actual performance gain of the proposed method, an implementation of both the conventional and the multilevel version of the Walk on Spheres algorithm was created.
By parallelizing the generation of samples, more thorough convergence results could be obtained.
Additionally, a version of the multilevel method that chooses the number of samples based on measurements was implemented.

In order to test a variety of different situations, multiple domains and boundary conditions were implemented.
The convergence rate of the error with respect to the work was measured, along with many other relevant quantities.
``Multilevel Walk on Spheres'' resulted in a reduction of the work by up to a factor 2. 

Using tighter bounds, it would perhaps be possible to find an upper bound for the variance convergence rate that better fits the measured behavior.
This would hopefully allow to analytically determine the number of samples per level.

It would be beneficial to generalize the ``Multilevel Walk on Spheres'' to other elliptic equations, for example the Poisson equation with nonconstant right-hand side.
This equation has applications in many fields, such as particle accelerator modeling or biochemical electrostatics, in which a Multilevel Monte Carlo formulation may be of use.
This would require the formulation of a process that contains multiple discretization levels \textit{along the path}, i.e., with level dependent maximal sphere radius. This would increase the expected work especially for fine discretizations and presumably increase the benefit of MLMC over MC. 

\bibliography{sources}
\bibliographystyle{abbrv}
\end{document}